\documentclass[12pt]{article}
\usepackage{latexsym,amsfonts,amssymb,mathrsfs,amsmath,amsthm,cite}
\usepackage[normalem]{ulem}

\usepackage{color}
\usepackage{colordvi,multicol}
\usepackage{geometry}
\newtheorem{thm}{\hskip\parindent\bf{Theorem}}
\newtheorem{cor}{\hskip\parindent\bf{Corollary}}
\newtheorem{lem}{\hskip\parindent\bf{Lemma}}

\newtheorem{rem}{\hskip\parindent\bf{Remark}}
\newtheorem{exa}{\hskip\parindent\bf{Example}}

\begin{document}
\title{\bf Lundberg-type inequalities
\\ for non-homogeneous
risk models
}

\renewcommand{\thefootnote}{}
\footnotetext{
CONTACT Junyi Guo    jyguo@nankai.edu.cn.   School of Mathematical Sciences,
 Nankai University,  Tianjin,  China.}

\author{\noindent{{Qianqian Zhou}$^{a}$, {Alexander Sakhanenko}$^{a,b}$   and  {Junyi Guo}$^{a}$}\\
\small {\noindent{$^a$School of Mathematical Sciences,
 Nankai University,  Tianjin,  China;}}\\
\small{ \noindent{$^b$Sobolev Institute of Mathematics, 4 Acad. Koptyug avenue,
 Novosibirsk, Russia}}}

\date{}
\maketitle

\vskip 0.5cm \noindent{\bf Abstract}\quad
In this paper, we investigate  the ruin probabilities  of  non-homogeneous risk models. By employing martingale method, the Lundberg-type inequalities of ruin probabilities of non-homogeneous renewal risk models are obtained under weak assumptions. In addition,  for the periodic and quasi-periodic risk models the adjustment coefficients of the Lundberg-type inequalities are obtained.
Finally, examples are presented to show that  estimations obtained in this paper are more accurate and the ruin probability in non-homogeneous risk models may be fast decreasing which is impossible for the  case of homogeneity.

\smallskip

\noindent {\bf Keywords:}\quad  Non-homogeneous risk model;  Martingale method; Ruin probability;  Lundberg-type inequality

\smallskip

\noindent {\bf Mathematics Subject Classification (2010)}\quad 91B30; 60K10

\section{Introduction and some basic results}
In classical risk theory, the probability of ruin  is a vital  index of  the robustness of an insurance company  and is also a useful tool for risk management. Therefore, it has been studied by many authors and is still attracting the  attention of many authors.
High probability of ruin  means that the insurance company is not stable, and then  the insurer should take corresponding measures to reduce  risks.
Thus the research  on ruin probability of risk model is a hot subject in risk theory.
The calculation of ruin probability is a classical problem in actuarial science. However, the exact value of ruin probability can only be calculated for exponential distribution or discrete distribution with finite values.
But   the upper bounds of run probabilities of risk models can be obtained. The insurer may employ the upper bound of  ruin probability to  evaluate the stability of  insurance company.

Let $R(t)$   be a risk reserve process, defined for all $t\ge 0,$ with non-random initial reserve $R(0) = u>0.$
 The ultimate ruin probability $\psi(u)$ of the risk model is defined as
\begin{gather}                                                        \label{aa1}
\psi(u):={\bf P}\Big[\inf_{t\geq 0}R(t)<0\Big| R(0)=u\Big].
\end{gather}
After the classic works of Lundberg \cite{LF03} and Cram$\acute{e}$r \cite{C30, C55}, many research works,   such as references
\cite{AA10, DC05, GH79, G91, RSST98},   have been devoted to studying   ruin probabilities of homogeneous risk models.

For classical  homogeneous risk models the most excellent result about the behavior of $\psi(u)$
 is the Lundberg inequality, which states that  under some appropriate assumptions (see  \cite{C30, C55} and \cite{LF03} for more details)
\begin{gather}                                                                                \label{aa9}
\psi(u)\leq e^{-Lu}\text{ for all } u\geq 0,
\end{gather}
The largest number $L$ in \eqref{aa9} is called the adjustment coefficient or Lundberg exponent.

In the  homogeneous risk models,  "identically distributed"  assumption is imposed both on  inter-occurrence times and claim sizes.  However,
both factors are influenced by the economic environment and the usual assumptions on claim sizes and inter-occurrence times may be too restrictive for practical use. Thus relaxing  the "identically distributed" assumption imposed  on claim sizes or (and) inter-occurrence times is necessary, which yields  non-homogeneous risk processes.
In Casta$\tilde{n}$er, et al. \cite{CCGLM13}, and  Lef$\grave{e}$vre  and Picard \cite{LP06},  the "identically distributed" assumption  imposed on claim sizes is relaxed, while in Bernackait$\dot{e}$ and $\breve{S}$iaulys \cite{BJ17},    Ignator  and Kaishev \cite{IK00} and  Tuncel and  Tank \cite{TT14}, the "identically distributed" assumption  imposed on inter-occurrence times is also relaxed.

Without the condition of identical  distribution on claim sizes and inter-occurrence times,
  there are  serious difficulties  in  evaluating  the probability of ruin. Hence, the related papers generally investigate the recursive formulas for the finite time ruin  probability under restrictive conditions. Such as in Bla$\check{z}$evi$\check{c}$ius, et al. \cite{BBE10},  the recursive formula of finite-time ruin probability of discrete time risk model with nonidentical distributed claims is obtained. In R$\check{a}$ducan, et al. \cite{RVZ15a}, the recursive formulae for the ruin probability at or before a certain claim arrival instant of  the non-homogeneous risk model are obtained. In this risk model, the claim sizes are independent non-homogeneous Erlang distributed  and independent of the inter-occurrence times, which are assumed to be i.i.d. random variables  following an Erlang or a mixture of exponential distribution. R$\check{a}$ducan, et al. \cite{RVZ15b} extended those recursive formulae to a more general case when the inter-occurrence times are i.i.d. nonnegative random variables following an arbitrary distribution.

Our aim is to investigate  non-homogeneous risk models and  obtain Lundberg inequality   for
the probability $\psi(u)$  under more general assumptions.

We  consider a class of risk reserve processes, which need not to be homogeneous, with the following properties.
\medskip

 \hspace*{0.5cm} (i) Process $R(t)$ may have positive jumps only at random or non-random times $T_1, T_2, \ldots$ such that
 \begin{gather*}
\forall k=1, 2, \ldots,\ T_{k+1}>T_{k} >T_0:=0 \text{ and } T_n\rightarrow\infty\ a.s.
\end{gather*}

\hspace*{0.5cm} (ii) Process $R(t)$ is monotone on each time interval $[T_{k-1}, T_k), \ k=1,2,\ldots,$ and
$R(0)=u>0.$

\medskip

{\bf Model A.} \  Assume that the \emph{kth}  claim $Z_k$ {occurs at}  time $T_k$, i.e.,
\begin{gather}                                                               \label{aa5}
-Z_k:= R(T_k)-R(T_k-0)\le 0,\ \theta_k:=T_k-T_{k-1}>0,\  k=1,2,\ldots.
\end{gather}
Suppose that on each interval $[T_{k-1}, T_k)$ {the premium  rate is} $p_k,$ i.e.,
\begin{gather}\label{aa6}
\forall t\in[T_{k-1}, T_k),\  R(t)-R(T_{k-1})=p_k(t-T_{k-1}), \  k=1,2,\ldots.
\end{gather}
Assume, also, that random vectors
\begin{gather*}\label{aa7}
(p_k, Z_k, \theta_k), \ k=1,2,\ldots
\end{gather*}
are mutually independent. Then conditions (i) and (ii) hold  and random variables
\begin{gather*}
Y_k=R(T_{k-1})-R(T_k)=Z_k-X_k=Z_k-p_k \theta_k,\ k=1,2,\ldots
\end{gather*}
are also mutually independent. Here we denote  by $X_k=R(T_k-0)-R(T_{k-1})=p_k \theta_k$ the total premium {collected by the insurer} over the time  interval $[T_{k-1}, T_k)$.

Here we call  model A is non-homogeneous renewal risk model and  it is clear that model A is more general than the classical compound Poisson and  renewal risk models  introduced by  Sparre  Andersen \cite{AE57}   in which the random vectors $(Z_k, \theta_k),\  k=1,2,\ldots$, are assumed to be i.i.d. and the premium  rate $p_k\equiv p>0$ is fixed, positive and non-random.

\begin{thm}\label{thm1}
Assume that conditions (i) and (ii) hold and the  random variables
\begin{gather}                                                         \label{aa3}
Y_k:=R(T_{k-1})- R(T_k),\ k=1,2,\ldots,
\end{gather}
 are mutually independent.
Then for any $ u>0$ and $h\ge 0,$ the ruin probability  $\psi(u)$   satisfies
\begin{gather}                                                                 \label{aa4}
\psi(u)=
{\bf P}[\sup_{k\geq 1}S_k>u]\leq e^{-hu}\sup_{k\geq1}{\bf E}e^{hS_k}
=e^{-hu}\sup_{k\geq1}\prod_{j=1}^k{\bf E}e^{hY_j},
\end{gather}
where $S_k=Y_1+Y_2+\ldots+Y_k.$

Moreover, the Lundberg inequality (\ref{aa9})  holds with $L=L(S_\bullet),$
where
\begin{gather}                                                                                        \label{aa11}
L(S_\bullet):=\sup\{h\ge0: \sup_{k\ge 1}{\bf E}e^{hS_k}\le 1 \}.
\end{gather}
\end{thm}

Note that   in   Theorem \ref{thm1}  the  generalization of  Lundberg inequality  for non-homogeneous  risk models  is obtained. We also have  the following  assertion  which follows immediately from Theorem \ref{thm1}.

\begin{cor}                                                              \label{CC1}
Under the  assumptions of Theorem \ref{thm1}, for any  $u>0,$ $\psi(u)$  satisfies that
\begin{gather}                                                                                        \label{aa12}
\psi(u)\leq \inf_{h\in[0, L(Y_\bullet)]}\{e^{-hu}{\bf E}e^{hY_1}\}\leq e^{-L(Y_\bullet) u}{\bf E}e^{L(Y_\bullet) Y_1}\leq e^{-L(Y_\bullet)u},
\end{gather}
with
${L(Y_\bullet)}:=\sup\{h\ge0:\sup_{j\geq 1}{\bf E}e^{hY_j}\leq 1\}\le L(S_\bullet).$
\end{cor}

 Homogeneous renewal risk model is a special case of the model from  Corollary \ref{CC1} when  random variables $Y_1, Y_2, \ldots$ are i.i.d..
In this case from \eqref{aa12} we have the classical  Lundberg inequality \eqref{aa9}
with  the adjustment coefficient or Lundberg  exponent $L=L(Y_1)$ given by
\begin{gather}                                                                                        \label{aa13}
L(Y_1):=\sup\{h\ge0: {\bf E}e^{hY_1}\leq 1\}.
\end{gather}
It can be seen that our inequality \eqref{aa12} is a little better than the Lundberg inequality in its classical form,
which can be found in  \cite{AA10, G91,RSST98},
 because
 in our variant~\eqref{aa12} of Lundberg inequality we do not exclude
the cases of
${\bf E}e^{L(Y_1)Y_1}<1$ and/or  ${\bf E}Y_1=-\infty.$

It follows from   \eqref{aa13} that ${\bf E}e^{hY_1}>1$ for all $h>L(Y_1)$. So, for i.i.d. random variables $Y_1, Y_2, \ldots$ we have that the right hand side in \eqref{aa4} is $+\infty$
 for all $h>L(Y_1)$.
Thus,
 for  classical renewal risk models  the adjustment coefficient $L(Y_1)$ is the
 natural boundary for possible values of the parameter $h$ in  inequality  \eqref{aa4} of  Theorem \ref{thm1}. But for  non-homogeneous  risk models
the situation may be significantly different because
 the optimal value  $h=h(u)$ of the parameter $h$ in  inequality  \eqref{aa4}  may be greater  than $L(S_\bullet)$ and may tends to $+\infty$ as $u\to\infty$. Moreover, in Examples \ref{exa3} and \ref{exa4} below we present random variables $Y_1, Y_2, \ldots$
corresponding  to risk models such that
\begin{gather}                                                        \label{aa15}
\psi(u)=o\big(e^{-Nu}\big) \text{  as  } u\to\infty \text{ for all }  N<\infty.
\end{gather}
 Moreover, in Example \ref{exa3} we have:
\begin{gather}                                                                  \label{aa16}
\forall u>0, \
\psi(2u)\leq
 e^{-u^{3/2}}.
\end{gather}
So, very fast decreasing of ruin probabilities is  possible in non-homogeneous cases.

Theorem \ref{thm1} allows us also to obtain a
generalization of the Lundberg inequality  that for any $u>0$
\begin{gather}                                                        \label{aa17}
 \psi(u)\leq Ce^{-Lu} \text{ with }\  C<\infty \text{ and } L>0.
\end{gather}
For example, for  the periodic risk model with period $l$ Theorem \ref{thm1} yields immediately that
\begin{cor}                                                              \label{CC3}
Suppose that there exists an integer $l\geq 1$  such that for all $n=1,2, \ldots$ random variables $Y_{n+l}$ and $Y_n$
are identically distributed. Then under assumptions   of Theorem \ref{thm1} inequality (\ref{aa4}) holds  with
\begin{gather*}
\sup_{k\geq 1}{\bf E}e^{hS_k}=
\max_{1\le k\le l}{\bf E}e^{hS_k}\text{ for each }  h\in[0,L(S_l)],
\end{gather*}
where
 $$L(S_l):=\sup\big\{h\ge 0:  {\bf E}e^{hS_{l}}\le 1     \big\}.$$

Moreover, for any $u>0$
$$\psi(u)\le \inf_{h\in [0, L(S_l)]}\{ e^{-hu}\max_{1\le k\le l}{\bf E}e^{hS_k}   \}\le C_1 e^{-L(S_l)u}$$
with $C_1=\max_{1\le k\le l}{\bf E}e^{L(S_l)S_k}.$

\end{cor}

There are a few works in which the Lundberg inequalities of ruin probabilities of non-homogeneous  models are studied.
In   Andrulyt$\dot{e}$, et al. \cite{AE15},  Kievinait$\dot{e}$  \& $\check{S}$iaulys \cite{KD18}, and Kizinevi$\check{c}$  \& $\check{S}$iaulys \cite{EJ18},  the Lundberg inequalities of non-homogeneous renewal risk models were obtained. In their models, they all assume  that claim sizes and  inter-occurrence times are both independent but not necessarily identically distributed.
We'll show in Remarks \ref{Rem3+} and \ref{Rem3} that their results are special cases of our Theorem \ref{thm1}.
For example, it is not possible for them to obtain inequalities with properties \eqref{aa15} or \eqref{aa16}. In the following, a general  non-homogeneous renewal risk model with interest rate is  studied. Thus
the results of the present  paper are complementary to the results in Andrulyt$\dot{e}$, et al. \cite{AE15},  Kievinait$\dot{e}$  \& $\check{S}$iaulys \cite{KD18}, and Kizinevi$\check{c}$  \& $\check{S}$iaulys \cite{EJ18}.

 The rest of the paper is organized as follows.
In Section \ref{section2}, we investigate the ruin probabilities of  non-homogeneous risk modes. The general risk model is studied then   the special case of periodic and quasi-periodic risk model is studied. Finally, some remarks of our results are presented. Examples are given in Section 3 to show that our estimations are more accurate and the probability of ruin   in  the non-homogeneous risk model may be fast decreasing which is impossible in homogeneous case.
Almost all proofs are gathered  in the last  section.

Later on we regularly use the fact  that expectations ${\bf E}e^{h S}\in(0, \infty]$ are everywhere defined for all random variables $S$ and all real number $h\ge 0$ but may take value of positive infinity.
By this reason  all inequalities of the form ${\bf P}[A]\le \text{const} \cdot {\bf E}e^{h S}$ make sense even we omit, for brevity, the assumption  ${\bf E}e^{h S}<\infty$.

Note also that for the probability $\psi(u)\le1$  the inequality $\psi(u)\le\psi^*(u)\le\infty$ means that $\psi(u)\le\min\{\psi^*(u),1\}\le1$.
Later on  we will use an agreement that
\begin{gather*}
{ E}+\text{const}=\infty   \text{ and }  \text{const}/{E}=0 \text{ when }   {E}=\infty.
\end{gather*}
All limits in this paper are taken with respect to $n\to\infty$ unless the contrary is specified.
And  we use random  variables $Y_1, Y_2, \ldots$
only  when they are mutually independent.

\section{Ruin probability for non-homogeneous risk models}\label{section2}

\noindent{\bf 2.1. General risk model }

Now we consider a more general class of risk reserve processes $R(t)$ which, together with properties (i) and (ii),  satisfy the following two assumptions.

\medskip
 \hspace*{0.5cm} (iii) For some non-random $r_1, r_2,\ldots$
 \begin{gather}\label{Riii}
\forall n\geq 1,\ R(T_k)\ge(1+\alpha_k)R(T_{k-1})-Y_k \text{ and } \alpha_k\ge r_k\ge0.
\end{gather}

\hspace*{0.5cm} (iv) Random variables $Y_k^*:=\frac{Y_k}{1+\alpha_k}, \ k=1,2,\ldots,$
are  mutually independent.
\medskip


\medskip

{\bf Model B.} \
 Instead of \eqref{aa6}
suppose that on each interval $[T_{k-1}, T_k)$ the surplus  accumulates    as follows
\begin{gather*}
R(t)=(1+\alpha_k(t))R(T_{k-1})+(1+\beta_k(t))p_k(t-T_{k-1}), \ k=1,2,\ldots,
\end{gather*}
where $(1+\alpha_k(t))R(T_{k-1})$ is the accumulated value of $R(T_{k-1})$ from $T_{k-1}$ to $t$ under interest rate  $\alpha_k(t)$
and $(1+\beta_k(t))p_k(t-T_{k-1})$ is the accumulated value of premiums collected from $T_{k-1}$ to $t$ under premiun rate $p_k$ and interest rate  $\beta_k(t)$.

 Assume again that claims $Z_k$ arrive to an insurer only at   time $T_k$, i.e.,~\eqref{aa5} holds.
Suppose  now that   processes $\alpha_k(t)\ge0$ and $\beta_k(t)\ge0$ are non-decreasing, and  random vectors
\begin{gather*}
(p_k, Z_k, \theta_k, \alpha_k(T_k), \beta_k(T_k) ), \ k=1,2,\ldots,
\end{gather*}
are mutually independent.
 Then conditions (i)-- (iv) hold  and random variables
\begin{gather}\label{A1}
Y_k=(1+\alpha_k(T_k))R(T_{k-1})-R(T_k)=Z_k-X(T_k)=Z_k -(1+\beta_k(T_k))p_k \theta_k,\ k=1,2,\ldots,
\end{gather}
are also mutually independent.
\medskip

 It is clear that the presented model is more general than model A.

For each $k=1, 2, \ldots$ introduce notations:
\begin{gather}\label{R+}
  v_k:=\prod_{j=1}^{k}\frac{1}{1+r_j}\text{ and }
S_k^*:=\sum_{j=1}^k v_{j-1}Y_j^*,
\end{gather}
 where $v_0:=1$.

\begin{thm}\label{thm2}
Under   assumptions (i)--(iv), for any $u>0$  and any $h\ge 0$ the following inequality
\begin{gather}
\label{R6+}
 \psi(u)\leq e^{-hu}\sup_{k\geq1}{\bf E}e^{hS_k^*}
=e^{-hu}\sup_{k\geq1}\prod_{j=1}^k{\bf E}e^{hv_{j-1}Y_j^*},
\end{gather}
 holds.

In addition, the Lundberg inequality (\ref{aa9}) holds with $L=L(S_\bullet^*)$
where
\begin{gather}\label{A3}
L(S_\bullet^*):=\sup\Big\{h\ge0: \sup_{k\ge 1}{\bf E}e^{hS_k^*}\le 1 \Big\}.
\end{gather}

\end{thm}

\begin{rem}
As an analogue of Corollary \ref{CC1} for   model A, it is easy to see that all the assertions of Corollary \ref{CC1} also hold for  model B with
$$ L(Y_\bullet)=\underline{L}:=\sup\Big\{h\ge0:\sup_{j\geq 1}{\bf E}e^{hv_{j-1}Y_j^*}\leq 1\big\}\le L(S_\bullet^*).$$
\end{rem}

  Theorem \ref{thm2} also yields  the following result.
\begin{cor}                                                              \label{C2+}
Suppose that for each $k\ge 1$
 \begin{gather*}
 Y_k=b_k\xi_k\text{ and }b_kv_{k-1}\le1,
\end{gather*}
where  random vectors $\{(\alpha_k,\xi_k),k=1,2,\ldots,\}$ are i.i.d..
In this case under assumptions (i)--(iii)
inequality \eqref{aa12} holds with
$$L(Y_\bullet)=L(Y_1^*):=\sup\{h\ge0: {\bf E}e^{hY_1^*}\leq 1\}.$$
In particular, for any $u>0$
\begin{gather}                                                                                        \label{c6+}
 \psi(u)\leq  e^{-\kappa u} \text{  if  }
{\bf E}e^{\kappa Y_1^*}=1.
\end{gather}
\end{cor}

Thus, when $b_n=(1+r)^{n}$ and $v_n=(1+r)^{-n}<1$ we obtain two generalizations \eqref{aa12} with $L(Y_\bullet)=L(Y_1^*)$ and \eqref{c6+} of the famous Lundberg inequality for the case when $|Y_n|=(1+r)^{n}|\xi_n|\to\infty$ almost surely and with high speed.

\begin{rem}                                                                                 \label{Rem4}
Two special cases of  inequality \eqref{c6+} are obtained in Corollaries 3.1 and 3.2 of Cai \cite{CJ02a}  when
$$
Y_k=X_k-Z_k \text{ or }
Y_k=(1+\alpha_k)X_k-Z_k
$$
with $r_k\equiv0$.
Earlier a simpler case with
$$
\alpha_k=r\ge0\text{ and }
Y_k=(1+r)X_k-Z_k,
$$
 is considered in  Yang \cite{YH99}.
Underline that the mutual independence of non-negative random variables $X_k$, $Z_k$ and $\alpha_k$ is essential for the proofs in  \cite{CJ02a, YH99}.

\end{rem}

\noindent{\bf 2.2. Periodic and quasi-periodic risk model }

Asmussen \& Rolski \cite{AR94} (see also Asmussen  \& Albrecher \cite{AA10}) studied a kind of risk process which happens  in a periodic environment. For the Lundberg-type inequality \eqref{aa17} they obtained that adjustment coefficient $L$ is the same as for the standard time-homogeneous Poisson risk process obtained by averaging  the parameters over a period.  In Corollary \ref{CC3} we have found the similar property for periodic risk models with period $l$  under assumptions of Theorem~\ref{thm1}. Now we present two more general results under conditions of Theorem \ref{thm2}.

\begin{cor}\label{corP7}
Suppose  that there exist integers   $l,m\geq 1$  and a real number  $L^*\geq 0$  such that for any $n\ge m$
\begin{gather}                                                                                                               \label{p1}
{\bf E}e^{L^*(S_{n+l}^*-S_{n}^*)}\le 1.
\end{gather}
Then under assumptions (i)--(iv), for any $h\in[0, L^*]$ inequality \eqref{R6+}   holds     with
\begin{gather}                                                                                                             \label{p2}
\sup_{k\geq 1}{\bf E}e^{hS_k^*}
= \max_{1\leq k\leq l+m-1}{\bf E}e^{hS_k^*}.
\end{gather}

Moreover, for any $u>0$
\begin{gather}\label{A4}
\psi(u)\le \inf_{h\in [0, L^*]}\{e^{-hu}\max_{1\leq k\leq l+m-1}{\bf E}e^{hS_k^*} \}\le C_2 e^{-L^*u}
\end{gather}
with $C_2=\max_{1\leq k\leq l+m-1}{\bf E}e^{L^*S_k^*}.$

\end{cor}

\begin{thm}                                                              \label{thm3}
Suppose that there exist an integer $l\geq 1$ and a real number $q_l>0$ such that for all $n=1,2, \ldots$ random variables $Y_{n+l}^*$ and $q_lY_n^*$
are identically distributed.
  Assume also that for each $n\ge 1$
\begin{gather}                                                                                                         \label{p5}
q_l v_l\le1 \text{ and }   r_{n+l}=r_n.
\end{gather}
And set
\begin{gather}\label{p6}
L(S_l^*):=\sup\big\{h\ge 0:  {\bf E}e^{hS_{l}^*}\le 1     \big\}.
\end{gather}
Then under assumptions (i)--(iv), for any $h\in[0,L(S_l^*)]$ inequality \eqref{R6+}   holds
with
\begin{gather}                                                                                                             \label{p7}
\sup_{k\geq 1}{\bf E}e^{hS_k^*}\le
\max_{0\leq k< l}{\bf E}e^{hS_k^*}.
\end{gather}
And  for any  $u>0$ the following Lundberg-type inequality holds
\begin{gather}\label{A5}
\psi(u)\le \inf_{h\in [0, L(S_l^*)]}\{e^{-hu}\max_{0\le k<l}{\bf E}e^{hS^{*}_k}   \}\le C_3 e^{-L(S_l^*)u}
\end{gather}
with $C_3=\max_{0\le k<l}{\bf E}e^{L(S_l^*)S^{*}_k}.$

Moreover, if $q_lv_l=1$ then \eqref{p2}
is also true with $m=1$ for all $h\in[0, L(S_l^*)]$.
\end{thm}

Note that the value in the right hand side of \eqref{p2}
may be less than 1. On the other hand,  the right hand side of \eqref{p7}
may not be less than 1 since ${\bf E}e^{hS_0^*}=e^0=1.$
Thus, Corollary \ref{corP7}
may give sharper estimates than Theorem \ref{thm3}
because condition \eqref{p1} may be stronger than the assumptions in Theorem \ref{thm3}.

The models satisfying conditions of Corollary \ref{corP7} or Theorem \ref{thm3} will be called quasi-periodic.
For us it is essential that for $l>1$ all such risk models are automatically non-homogeneous.
Models satisfying assumptions of Theorem \ref{thm3} with $q_lv_l=1$ are naturally  called  periodic.
Periodic model from Corollary \ref{CC3} is a  special case  of Theorem \ref{thm3}  with $q_l=v_l=1=m$ when $S_k^*=S_k$ and  \eqref{p2} is true.

\noindent{\bf 2.3. General remarks. }

\begin{rem}                                                                                 \label{Rem5}
Theorem \ref{thm1} is an evident special case of Theorem  \ref{thm2} when
$\alpha_k=r_k=0$ for all $k=1,2,\dots$. Hence, all corollaries from Theorem  \ref{thm2}, which are presented below, may be considered also as corollaries from Theorem \ref{thm1}.
\end{rem}

\begin{rem}                                                                                 \label{Rem5}
Theorems \ref{thm1} and \ref{thm2} and all their corollaries remains valid also for non-homogeneous discrete-time risk models  when claims arrive {at} non-random times $T_1,T_2,\dots$.
\end{rem}

\begin{rem}                                                                                 \label{Rem3+}
In  \cite{AE15, KD18},     the authors use the trivial inequality
\begin{gather*}
\forall u>0, \  \forall h\ge0, \   {\bf P}[\sup_{k\geq 1}S_k>u]\leq\sum_{k=1}^\infty
{\bf P}[S_k>u]\leq e^{-hu}\sum_{k=1}^\infty {\bf E}e^{h S_k}
\end{gather*}
{instead of} our  sharper estimate \eqref{aa4}. It is clear that the results in  \cite{AE15, KD18}  will be improved automatically by using our  estimate \eqref{aa4}.
\end{rem}

\begin{rem}                                                                                 \label{Rem3}
In the proof of Theorem 4 in   \cite{EJ18} it is shown (under several additional assumptions)
that for any $u>0$
\begin{gather}                                                                                        \label{c1-}
\psi(u)\leq \inf_{h\in[0, L(Y_\bullet)]}\{e^{-hu}\sup_{i\geq 1}{\bf E}e^{hY_i}\}.
\end{gather}
It is clear that the inequality \eqref{aa12} is better than \eqref{c1-}.
It also follows from Example \ref{exa3} below that our estimate \eqref{aa4} from Theorem \ref{thm1}
may be qualitatively  better than the estimate \eqref{c1-}.
\end{rem}

\begin{rem}                                                                                 \label{Rem0}
Note that all independent random variables $\{Y_n\}$ (with all possible distributions) may appear
in Theorem \ref{thm1} since we everywhere may consider model A  with values
\begin{gather*}
Z_n=Y_n^+:=\max\{Y_n,0\}, \ X_n=p_n\theta_n=Y_n^-:=\max\{-Y_n,0\},\ p_n=1.
\end{gather*}
 For this reason, in Examples \ref{exa0}--\ref{exa4}   below we do not construct risk models but only introduce random variables $\{Y_n\}$ that appear in risk models.

Similarly, all random variables $\{Y_n^*\}$ (with all possible distributions)
may appear in Theorem \ref{thm2}
with all possible {real numbers $r_n\ge 0$}. Indeed, we may  consider model  B  with
\begin{gather*}
Z_n=(1+r_n)\max\{Y_n^*,0\}, \ X_n(T_n)=p_n\theta_n=\max\{-Y_n^*,0\},\ p_n=1,
\end{gather*}
and with $ \alpha_n(T_n)= \beta_n(T_n) =r_n\ge0$ for all $n=1,2,\ldots$.
\end{rem}

\section{Examples}
In this section,  special examples are presented to show that our estimations are more accurate and the ruin  probability in non-homogeneous risk model may be fast decreasing which is impossible in homogeneous case.

\begin{exa}\label{exa0}
Let $Y_1, Y_2, \ldots$  be independent normal
random variables with
\begin{gather}                                                                                                         \label{p11}
Y_n \sim N(a_n, 1)   \text{ and }  a_n+a_{n+1}\le a_1+a_2=-1,\     n=1, 2, \ldots.
\end{gather}
It is easy to calculate that
\begin{gather}                                                                                                         \label{p12}
{\bf E}e^{hY_n}=e^{ha_n+\frac{h^2}{2}} \text{ and } {\bf E}e^{hS_n}={\bf E}e^{h\sum_{i=1}^na_n+\frac{n}{2}h^2},\ n=1,2, \ldots.
\end{gather}
Hence,
$$L(S_2)=\sup\{h\ge 0: {\bf E}e^{hS_2}\le 1  \}=\sup\{h\ge 0:e^{-h+h^2}\le1 \}=1.$$
Thus, we have from \eqref{p11} and \eqref{p12} with $h=1$ that
\begin{gather}                                                                                                               \label{p1+}
{\bf E}e^{S_{n+2}-S_{n}}=e^{a_{n+1}+a_{n+2}+1}\le e^0= 1 ,\ \forall  n\ge 0.
\end{gather}

Comparing \eqref{p1} and \eqref{p1+} we obtain that random variables $\{Y_n\}$ from \eqref{p11}
satisfy  all conditions of  Corollary \ref{corP7}
with $l=2$, $m=1$ and $L^*=1$.
So, we have from \eqref{R6+} and \eqref{p2} with $h=L^*=1$ that for any $u>0$
\begin{gather}                                                                                                         \label{p14}
 \psi(u)\leq \max \{{\bf E}e^{S_1},  {\bf E}e^{ S_2}\}e^{-u}=e^{(a_1+1/2)^+-u}, 
\end{gather}
because
$$\max \{{\bf E}e^{S_1},  {\bf E}e^{ S_2}\}=\max \{e^{a_1+1/2},   1\}=e^{\max \{a_1+1/2,   0\}}.$$

Inequality \eqref{p14} allows us to make several conclusions about  risk models
with random variables $\{Y_n\}$ from \eqref{p11}. First, if $a_1\leq -1/2,$   we can obtain the Lundberg  inequality~(\ref{aa9})  with $L=L^*=1.$ Second, if  $a_1>-1/2,$  we can prove the generalization  (\ref{aa17}) of the Lundberg  inequality  with $C=e^{a_1+1/2}>1.$
Third, in the case of $a_1\ge0$  we have
\begin{gather*}
\sup_{i\geq 1}{\bf E}e^{hY_i}\ge {\bf E}e^{hY_1}=e^{a_1h+h^2/2}\ge e^{h^2/2}\ge e^0=1,\ \forall h\ge 0.
\end{gather*}
Thus,  it follows from \eqref{aa11} that in this case $L(Y_\bullet)=0$ and, hence, inequality \eqref{c1-} allows us to obtain only trivial estimate $\psi(u)\le1$.
So, (see Remark \ref{Rem3}) Theorem 4 from \cite{EJ18} does not work in this case, whereas, our results yield the estimate
(\ref{aa17})  with $C=e^{a_1+1/2}<\infty.$
\end{exa}

\begin{exa}(see \cite[Example 1]{KD18})                                                               \label{exa1}
Suppose   $Y_1, Y_2, \ldots$  are independent random  variables  such that:\\
\hspace*{0.5cm} $\bullet$ $Y_i$ are uniformly distributed on interval $[0,2]$  for  $i\equiv 1  \mod 3;$\\
\hspace*{0.5cm}  $\bullet$  $Y_i$  are uniformly distributed on interval $[-1,0]$  if $i\equiv 2  \mod 3;$\\
\hspace*{0.5cm}  $\bullet$  $\overline{F}_{Y_i}(x)=1_{(\infty, -2)}(x)+e^{-x-2}1_{[-2, \infty)}(x)$  when $i\equiv 0  \mod 3.$

It is a periodic risk model with $l=3$ and it is easy to see  that this model satisfies all the conditions of Corollary \ref{CC3}.  We can calculate that
for all $0<h<1$
\begin{gather*}
{\bf E}e^{hY_2}=\frac{1-e^{-2h}}{2h}=e^{-2h}{\bf E}e^{hY_1}<1,\
{\bf E}e^{hY_3}= \frac{e^{-2h}}{1-h}.
\end{gather*}
Then for $h_0=2/3$
\begin{gather*}
{\bf E}e^{h_0S_3}=\frac{(1-e^{-2h_0})^2}{(2h_0)^2(1-h_0)}<1,
\ {\bf E}e^{h_0Y_1}=\frac{e^{2h_0}-1}{2h_0}<2.2.
\end{gather*}
  Thus, we have from   Corollary \ref{CC3} that
$$
\psi(u)\leq  \max \{{\bf E}e^{S_1},  {\bf E}e^{ S_2},  {\bf E}e^{ S_3}\}e^{-h_0u}<2.2e^{-\frac{2}{3}u}.
$$
So, we have proved that
\begin{gather}                                                                       \label{p20+}
\psi(u)\leq \psi^{*}_1(u):=\max \{1, 2.2e^{-\frac{2}{3}u}\}
\leq  e^{(1-\frac{2}{3}u)^-},\ \forall u\geq 0.
\end{gather}

Remind that in  \cite{KD18}  the following bound is obtained
\begin{gather}                                                                              \label{p20}
\psi(u)\leq \psi_1^{\star}(u):=\min \big\{1, 1502e^{-0.01269u}   \big\},\ \forall u\geq 0.
\end{gather}
It is clear that estimate \eqref{p20+} is more accurate than  \eqref{p20}. For example,
$$\psi^{*}_1(576)< e^{-381}<10^{-165} \text{ whereas } \psi^{\star}_1(576)=1.$$

\end{exa}


Here we are going to present examples  of  risk models with property (\ref{aa15})
which is impossible in homogeneous case.

\begin{exa}                                                                \label{exa3}
Consider again
independent random variables $Y_1, Y_2, \ldots$ with normal distributions from Example \ref{exa0}.
First, suppose that they are i.i.d. with $ N(-1/2,1)$.
In this case condition (\ref{p11}) holds and $L(Y_1):=\sup\{h\ge0: {\bf E}e^{hY_1}\leq 1\}=\sup\{h\ge0:e^{-h/2+h^2/2}\leq 1\}=1$.
Thus we have the Lundberg inequality (\ref{aa9}) with $L=1$ and it is the upper boundary for values $h$ which we use in inequality~(\ref{aa4}).

Suppose now that condition (\ref{p11}) takes the form
\begin{gather}                                                                                                         \label{p11+}
Y_n \sim N(a_n, 1) \text{ with }  a_n=(1-2n)/4,\    n=1, 2, \ldots.
\end{gather}
In this case $\sum_{i=1}^na_i=-n^2/4$ and we obtain from (\ref{p12}) that for  any $h\ge 0$  and each $n=1,2, \ldots$
\begin{gather*}                                                                                                         \label{p12+}
 {\bf E}e^{hS_n}=e^{-hn^2/4+nh^2/2}=e^{h^3/4-h(n-h)^2/4}\le e^{h^3/4}.
\end{gather*}
Hence, we have from  (\ref{aa4}) that for any $u>0$  and $h\ge 0$
\begin{gather}\label{p23}
\psi(u)\leq  e^{-hu}\sup_{n\geq1}{\bf E}e^{hS_n}\le e^{-hu+h^3/4}.
\end{gather}
With $h=2\sqrt{u/3}$ we obtain from (\ref{p23}) that for any $u>0$
\begin{gather*}                                                                  \label{p30}
\psi(u)\leq  e^{-4(u/3)^{3/2}}=  e^{-cu^{3/2}}  \text{ where } c^2=16/27.
\end{gather*}
\end{exa}
So, we obtain an example  of a risk model with property (\ref{aa16}) mentioned in the introduction.

\begin{exa} (see \cite[Example 2]{KD18})                                                               \label{exa4}
Suppose that $Y_1, Y_2, \ldots$  are independent random variables  with
$$
{\bf P}[Y_n=1]=\frac{1}{n+1}, \ {\bf P}[Y_n=-1]=1-{\bf P}[Y_n=1],\  n=1,2, \ldots.$$
In this case for all $h$
$$
{\bf E}e^{h Y_n}=e^{h}\frac{1}{n+1}+e^{-h}\Big(1-\frac{1}{n+1}\Big)=1+\frac{(1-e^{-h})(e^{h}-n)}{n+1}.
$$
So, ${\bf E}e^{h Y_n}\le1$ if and only if $n\ge e^{h}$. Hence, for each $h>0$
\begin{gather}\label{p25}
\sup_{n\geq1}{\bf E}e^{hS_n}={\bf E}e^{hS_m}  \text{ iff }  m+1\ge e^{h}\ge m.
\end{gather}
Thus, with $m$ from \eqref{p25}
$$
\sup_{n\geq1}{\bf E}e^{hS_n}={\bf E}e^{hS_m} <(e^{h})^m\le e^{he^{h}},
$$
and, using \eqref{aa4}, we obtain
\begin{gather}\label{p27}
\forall h,\ u>0, \
\psi(u)\leq  e^{-hu}\sup_{n\geq1}{\bf E}e^{hS_n}\le e^{-hu+he^{h}}.
\end{gather}
With $h=\log(u/2)>0$ we have from (\ref{p27}) that
\begin{gather}                                                                  \label{p28}
\forall u>2,  \
\psi(u)\leq \psi^{*}_2(u):= e^{-(u/2)\log(u/2)}=  \Big(\frac2u\Big)^{u/2}.
\end{gather}
So, we obtain  another example  of a risk model with property (\ref{aa15}).

Remind that in   \cite{KD18} the following bound is given
\begin{gather}\label{p29}
\forall u\geq 0,\ \psi(u)\leq \psi^{\star}_2(u):= \min\big\{1, 178e^{-\frac{u}{20}}    \big\}.
\end{gather}
It is clear that estimate \eqref{p28} is more accurate than  \eqref{p29}. For example,
$$\psi^{*}_2(103)< 10^{-88}\text{ in \eqref{p28}, whereas } \psi_2^{\star}(103)=1 \text{ in \eqref{p29} }.$$

\end{exa}

\section{Proofs }
\subsection{Key Lemma}
Before proving  our  main  results, we first  introduce  the following key lemma.
\begin{lem}\label{L1}
If  random  variables $X_1, X_2, \ldots$  are mutually independent, then for any real number $w$, any $h\ge 0$  and any $n\ge 1$
\begin{gather}                                                                        \label{m1}
 {\bf P}[\max_{1\leq k\leq n}W_k>w]\leq e^{-hw}\max_{1\leq k\leq n}{\bf E}e^{hW_k},
\end{gather}
where $W_k=X_1+X_2+\ldots+X_k.$  In addition, for any $w$  and any $h\ge 0$
\begin{gather}                                                                       \label{m2}
  {\bf P}[\sup_{k\geq 1}W_k>w]\leq e^{-hw}\sup_{k\geq 1}{\bf E}e^{h W_k}.
\end{gather}
\end{lem}

{\bf Proof.}
If $M_n(h):=\max_{1\leq k\leq n}{\bf E}e^{hW_k}=\infty$ then the inequality \eqref{m1} is obvious.
So, suppose that $0<M_n(h)<\infty$ and note that in this case the following sequence
\begin{gather*}
\mu_k=\frac{e^{hW_k}}{{\bf E}e^{hW_k}},\   k=0,1,2,\ldots, n, \text{ with }\mu_0=1,
\end{gather*}
is a martingale. This fact is  evident, since for all $k\ge0$
\begin{gather*}
\mu_k=\mu_{k-1}\frac{e^{hX_k}}{{\bf E}e^{hX_k}}  \text{ and }
{\bf E}\Big[ \frac{e^{hX_k}}{{\bf E}[e^{X_k}]} \Big| \mu_1, \ldots, \mu_{k-1}\Big]=1.
\end{gather*}

Thus,  by maximal  inequality  for  martingale
\begin{eqnarray*}
\forall x>0,\ {\bf P}\Big[\max_{1\leq k\leq n}\frac{e^{hW_k}}{{\bf E}e^{hW_k}}>x\Big]
=  {\bf P}\Big[\max_{1\leq k\leq n}\mu_k>x\Big]\leq \frac{{\bf E}\mu_n}x=\frac1x.
\end{eqnarray*}
Hence, with $x=e^{hw}/M_n(h)$ we obtain that
\begin{eqnarray*}
{\bf P}\Big[\max_{1\leq k\leq n}W_k>w\Big]
&=&{\bf P}\Big[\max_{1\leq k\leq n}\frac{e^{hW_k}}{M_n(h)}> x=\frac{e^{hw}}{M_n(h)}\Big]\nonumber\\
&\leq& {\bf P}\Big[\max_{1\leq k\leq n}\frac{e^{hW_k}}{{\bf E}e^{hW_k}}>x\Big]
\leq \frac1x=e^{-hw}M_n(h).
\end{eqnarray*}
So, inequality \eqref{m1} is proved.

Note that
$\max_{1\leq k\leq n}W_n\uparrow \sup_{k\geq 1} W_k$. Hence
$${\bf P}\Big[\sup_{k\geq 1}W_k>w\Big]
=\lim_{n\rightarrow\infty}{\bf P}\Big[\max_{1\leq k\leq n}W_k>w\Big]
\leq e^{-hw}\sup_{n\ge1}M_n(h)
= e^{-hw}\sup_{k\geq 1}{\bf E}e^{h W_k}.$$
Thus, Lemma \ref{L1} is proved.

\subsection{ Proof of Theorem \ref{thm1}}
It follows from (i) and (ii) that
\begin{gather}                                                                             \label{Ri}
\inf_{t\ge0}R(t)=\inf_{k\ge1}\inf_{t\in [T_{k-1}, T_k]}R(t)
=\inf_{k\ge1}\min\{R(T_{k-1}), R(T_k)\}=\inf_{k\ge0} R(T_k).
\end{gather}
Next,  from  the definition of $Y_k$ in  \eqref{aa3}  and using telescoping sum  we have    that
$$
R(T_k)= R(T_0)-\sum_{j=1}^kY_j=u-S_k,\ k=0,1,2,\ldots,
$$
with $S_k=Y_1+Y_2+\ldots +Y_k$ and $S_0=0$. Thus,
$\inf_{t\ge0}R(t)= u-\sup_{k\ge0}S_k$,
and hence, by   the definition of $\psi(u)$ in  \eqref{aa1} for any $u>0$  the following equality holds,
\begin{gather}                                                                             \label{m12}
 \psi(u)={\bf P}[u-\sup_{k\ge0}S_k<0]={\bf P}[\sup_{k\ge0}S_k>u].
\end{gather}

Since $S_0=0$, we have from \eqref{m12} that for any $u>0$
\begin{gather}                                                                             \label{m13}
 \psi(u)={\bf P}[\sup_{k\ge1}S_k>u].
\end{gather}
From \eqref{m13} and estimate \eqref{m2}  of Lemma \ref{L1} we obtain  the first inequality of \eqref{aa4}. And since $Y_1, Y_2, \ldots$  are independent, then all the assertions of (\ref{aa4}) are proved.

In addition, from (\ref{aa4}) and the definition of $L(S_\bullet)$  in (\ref{aa11}) we have that
\begin{eqnarray*}
\psi(u)&\le& \inf_{h\ge 0}\{e^{-hu}\sup_{k\ge 1}{\bf E}e^{hS_k}   \}\\
&\le& e^{-L(S_\bullet)u}\sup_{k\ge 1}{\bf E}e^{L(S_\bullet) S_k}\\
&\le& e^{-L(S_\bullet)u},
\end{eqnarray*}
i.e., the Lundberg inequality (\ref{aa9}) holds with $L=-L(S_\bullet).$

Thus all the assertions are proved.

\subsection{ Proof of Theorem \ref{thm2}.}
For each $k=1, 2, \ldots,$ introduce notations:
\begin{gather}                                                       \label{c11}
  v_k^*:=\prod_{j=1}^{k}\frac{1}{1+\alpha_j}\text{ and }
S_k^{**}:=\sum_{j=1}^k v_{j}^*Y_j=\sum_{j=1}^k v_{j-1}^*Y_j^*,
\end{gather}
with $v_0^*:=1$.

\begin{lem}                                                  \label{L-2}
If  $v_0=v_0^*=1$, $S_0^{**}=S_0^*=0$ and $\alpha_k\ge r_k\ge0$ for all $k\ge0$, then for any $n\ge 0$
\begin{gather}                                                                  \label{c12}
\max_{0\le k\le n}S_k^{**}\le \max_{0\le k\le n}S_k^*.
\end{gather}
\end{lem}
{\bf  Proof.}
Since $\alpha_j\geq r_j\geq 0$ for $j=1,2, \ldots,$  we have that
\begin{gather}                                                                  \label{c13}
\forall j\geq 1,\
c_{j}:=\prod_{i=1}^{j}\frac{1+r_i}{1+\alpha_i}
=c_{j-1}\frac{1+r_j}{1+\alpha_j}\ge c_{j-1}.
\end{gather}
Thus, real numbers $\{c_j\}$ have the following property
\begin{gather}                                                                  \label{c14}
\forall j\geq 1,\
1=c_0\ge c_1\ge \ldots\ge c_{j-1} \ge c_j>0.
\end{gather}

Next, from \eqref{c11} and \eqref{c13}  we have for any  $j=1, 2, \ldots$ that
\begin{gather}                                                                  \label{c15}
Y_j^* v^{*}_{j-1}
=Y^{*}_j \prod_{i=1}^{j-1} \frac{1}{1+r_i}  \prod_{i=1}^{j-1}\frac{1+r_i}{1+\alpha_i}=Y^{*}_j v_{j-1}  c_{j-1}=(S^{*}_j -S^{*}_{j-1})  c_{j-1}.
\end{gather}
Now, substituting  \eqref{c15} into  \eqref{c11} we obtain that for all $k\le n$
\begin{eqnarray*}
S^{**}_k&=&
\sum_{j=1}^k v_{j-1}^*Y_j^*=\sum_{j=1}^k
(S^{*}_j -S^{*}_{j-1})  c_{j-1}
=c_{k-1} S^{*}_k+\sum_{j=1}^{k-1}(c_{j-1}-c_j)S^{*}_j\nonumber\\
&\leq& c_{k-1} M_n^*+\sum_{j=1}^{k-1}(c_{j-1}-c_j)M_n^*=c_0 M_n^*=M_n^*:=\max_{1\leq j\leq n}S_j^*,
\end{eqnarray*}
where
we also use \eqref{c14}.
So, inequality \eqref{c12} is proved.

\medskip
{\bf  Proof of Theorem \ref{thm2}.}
Multiplying \eqref{Riii} by $v_k^*$ we obtain for any $k\ge 1$ that
 \begin{gather*}
 v_k^*R(T_k)\ge v_{k-1}^*R(T_{k-1})-v_k^*Y_k=v_{k-1}^*R(T_{k-1})-v_{k-1}^{*}Y_k^{*}.
\end{gather*}
Hence, by induction for any $k\ge 1$
\begin{gather*}
v_k^*R(T_k)\ge v_0^* R(T_0)-\sum_{j=1}^k v_j^*Y_j=u-S_k^{**}.
\end{gather*}
 This fact and Lemma \ref{L-2} imply  for any $n\ge 0$ that
\begin{gather}                                                                  \label{c18}
\min_{0\le k\le n}v_k^*R(T_k) \ge u-\max_{0\le k\le n}S_k^{**}\ge u-\max_{0\le k\le n}S_k^*.
\end{gather}

Inequality \eqref{c18} with $v_k^*>0$ implies  that for any $u>0$ and each $n=0,1,2,\ldots$
\begin{gather*}
 {\bf P}[ \min_{0\le k\le n}R(T_k)<0]={\bf P}[ \min_{0\le k\le n}v_k^*R(T_k)<0]
\le{\bf P}[u-\max_{0\le k\le n}S_k^*<0]={\bf P}[\max_{0\le k\le n}S_k^*>u].
\end{gather*}
Since $S^{*}_0=0$ and $\{S_k^*\}$ are
sums of independent random variables $\{v_{j-1}Y^{*}_j\},$ from Lemma  \ref{L1} for any  $u>0$  and each  $ n=1,2,\ldots$ we have that
\begin{gather}\label{A2}
{\bf P}[ \min_{0\le k\le n}R(T_k)<0]\le{\bf P}[\max_{1\le k\le n}S_k^*>u]\le e^{-hu}\max_{1\le k\le n}{\bf E}e^{hS^{*}_k}.
\end{gather}

On the other hand, equality \eqref{Ri} again follows from (i) and (ii)  for all $n\ge0$.
Then  taking limit on both sides of \eqref{A2}  as $n\to\infty$ we obtain the first inequality  in  \eqref{R6+}.
The equality in \eqref{R6+} directly comes from the definition of $S^{*}_k$ in (\ref{R+})  and the independence of $Y^*_j.$

The rest of proof is similar to Theorem \ref{thm1}.

\subsection{ Proof of Corollary \ref{corP7}}
For all  $n\ge m+l,$  since $Y^*_1, Y^*_2, \ldots$ are independent then the  random variables  $\Delta^*_{n,l}:=S_{n}^* -S_{n-l}^*$  and $S_{n-l}^*$  are independent. Hence, by Jensen's inequality for any $h\in [0, L^*]$  we    have from \eqref{p1} that
\begin{gather}                                                                  \label{c21}
{\bf E}e^{h\Delta^*_{n,l}}\le \Big({\bf E}e^{L^*\Delta^*_{n,l}}\Big)^{h/L^*}\le1\text{ and }
{\bf E}e^{hS^*_n}={\bf E}e^{h\Delta^*_{n,l}}\cdot {\bf E}e^{hS^*_{n-l}}\le {\bf E}e^{hS^*_{n-l}}.
\end{gather}
Since \eqref{p1} holds for any $n\ge m,$  then we can do (\ref{c21}) again  that
$${\bf E}e^{hS_{n-l}^*}\le {\bf E}e^{hS^*_{n-2l}}$$
if $n-l\ge m+l,$
otherwise ${\bf E}e^{hS^*_n}\le \max_{1\le k\le m+l-1}{\bf E}e^{hS^*_k}.$

Using induction with respect to $n$ it is not difficult   to see that we can do (\ref{c21}) for $i$ times until $n-il<m+l,$ i.e., $n-il\le m+l-1,$  with
$${\bf E}e^{hS^*_n}\le {\bf E}e^{h S^*_{n-il}}.$$

Then we obtain  that for any $n\ge m+l$
\begin{gather*}
 {\bf E}e^{hS^*_n}\le \max_{1\le k\leq m+l-1}{\bf E}e^{hS^*_k}.
\end{gather*}
Hence
$$\max_{1\le k\leq m+l-1}{\bf E}e^{hS^*_k}\le\sup_{n\ge 1}{\bf E}e^{hS^*_n}\le \max_{1\le k\leq m+l-1}{\bf E}e^{hS^*_k}.$$
So, equality \eqref{p2} is proved.

From (\ref{p2}) and the results in Theorem \ref{thm2}  we see that for any $h\in [0, L^*]$ and $u>0$
\begin{eqnarray*}
\psi(u)&\le& e^{-hu}\max_{1\le k\le l+m-1}{\bf E}e^{hS^*_k}\\
&\le& \inf_{h\in [0, L^*]}\{ e^{-hu}\max_{1\le k\le l+m-1}{\bf E}e^{hS^*_k}  \}\\
&\le& e^{-L^*u}\max_{1\le k\le l+m-1}{\bf E}e^{L^* S^*_k}.
\end{eqnarray*}
Therefore, the result in (\ref{A4})  is proved with $C_2=\max_{1\leq k\leq l+m-1}{\bf E}e^{L^*S_k^*}.$

\subsection{ Proof of Theorem \ref{thm3}.}
It is clear that each integer $n\ge1$ may be represented in the following way
\begin{gather}                                                                  \label{c24}
 1 \le n=il+k \text{ for some } i\ge0 \text{ and some } 1\le k\le l.
\end{gather}

\begin{lem}\label{lem4}
Under assumptions of Theorem \ref{thm3}
for all  $i=0,1,2, \ldots$ and $k=1,2, \ldots$, the following random variable
\begin{gather}                                                                                                         \label{c26}
\Delta_{i,k}:=S^*_{il+k}-S^*_{il}=\sum_{j=1}^{k}v_{ik+j-1}Y_{ik+j}^*
\end{gather}
is identically distributed with $(q_lv_l)^iS^*_{k}$.
In particular,  for any $ h\in[0,L(S_l^*)]$
\begin{gather}                                                                                                             \label{c27}
\sup_{n\geq 1}{\bf E}e^{hS_n^*}\le
\max_{1\leq k\leq l}\sup_{i\geq 0}{\bf E}e^{h(q_lv_l)^iS_k^*},
\end{gather}
where $L(S_l^*)$ is defined in (\ref{p6})
\end{lem}

{\bf Proof. }
It follows from  condition \eqref{p5} that
\begin{gather*}
v_{j-1}=\prod_{i=1}^{j-1}\frac{1}{1+r_i}
=\prod_{i=1}^{l}\frac{1}{1+r_i}\cdot\prod_{i=l+1}^{j-1}\frac{1}{1+r_i}
=v_{l}\cdot\prod_{i=1}^{j-l-1}\frac{1}{1+r_i}=v_{l}v_{j-l-1}.
\end{gather*}
Hence,  for all $n=1,2, \ldots$, random variables $v_{n+l-1}Y_{n+l}^*$ and $(q_lv_l)v_{n-1}Y_n^*$
are identically distributed because $Y_{n+l}^*$ and $q_lY_n^*$
are identically distributed by the assumption of Theorem~\ref{thm3}.
For $n=il+k,$ using induction with respect to $i$ it is not difficult to obtain that
$v_{il+j-1}Y_{il+j}^*$ and $(q_lv_l)^iv_{j-1}Y_j^*$ are identically distributed
for all  $i=0,1,2, \ldots$ and $j=1,2, \ldots$.
As a result,  $\sum_{j=1}^k v_{il+j-1}Y^*_{il+j}$ are identically distributed with $\sum_{j=1}^k (q_l v_l)^iv_{j-1}Y^*_j=(q_l v_l)^i\sum_{j=1}^k v_{j-1}Y^*_k=(q_l v_l)^iS^*_k.$

So, using definition (\ref{c26}), we have the first assertion of the lemma that
$\Delta_{i,k}$ and $(q_lv_l)^iS^*_{k}$  are identically distributed.

Using telescoping sum, we  note that $S^*_{il+k}=\sum_{m=0}^{i-1}\Delta_{m,l}+\Delta_{i,k}$. Hence
\begin{gather}                                                                                                             \label{c29}
{\bf E}e^{hS^*_{il+k}}=\prod_{m=0}^{i-1}{\bf E}e^{h\Delta_{m,l}}\cdot {\bf E}e^{h\Delta_{i,k}}
=\prod_{m=0}^{i-1}{\bf E}e^{h(q_lv_l)^mS^*_{l}}\cdot {\bf E}e^{h(q_lv_l)^iS^*_{k}}.
\end{gather}
Since $0\le h(q_lv_l)^m\le h\le L(S_l^*)<\infty$,  for any $h\in[0,L(S_l^*)]$ we have that
\begin{gather*}                                                                                                             \label{c29+}
{\bf E}e^{h(q_lv_l)^mS^*_{l}}
\le \big({\bf E}e^{L(S_l^*)S^*_{l}}\big)^{h(q_lv_l)^m/L(S_l^*)}\le1
\end{gather*}
by definition of $L(S_l^*)$.
From this inequality and  \eqref{c29}, for any $ h\in[0,L(S_l^*)]$ and each $i\ge 0$  and $k\ge 1$ we obtain
\begin{gather*}                                                                                                             \label{c29++}
{\bf E}e^{hS^*_{il+k}}\le{\bf E}e^{h(q_lv_l)^iS^*_{k}}.
\end{gather*}
The latter yields (\ref{c27}) if only we remind that
 each integer $n\ge1$ may be represented in the form (\ref{c24}).

Thus the lemma  is proved.

\medskip
{\bf Proof of Theorem \ref{thm3}.}
If $q_lv_l=1,$ then  for any $ h\in[0,L(S_l^*)]$  we have from  \eqref{c27} that
$$
\sup_{k\geq 1}{\bf E}e^{hS_k^*}\le\max_{1\leq k\leq l}{\bf E}e^{hS_k^*}.
$$
So, the second assertion of  Theorem \ref{thm3} is proved.

To prove the first one note that for $(q_lv_l)^i\le1$ and $1\le k\le l$ we have
\begin{gather}                                                                                                             \label{c30}
{\bf E}e^{h(q_lv_l)^iS_k^*}\le \big({\bf E}e^{hS_k^*}\big)^{(q_lv_l)^i}
\le \max\{1,{\bf E}e^{hS_k^*}\}\le
\max_{0\leq k\leq l}{\bf E}e^{hS_k^*} =\max_{0\leq k< l}{\bf E}e^{hS_k^*} ,
\end{gather}
because $S_0^*=0$ and ${\bf E}e^{hS_0^*}=1$.
Inequality \eqref{c30} together with \eqref{c27} implies \eqref{p7} under assumptions of Theorem \ref{thm3}.

The rest of proof is similar to Corollary \ref{corP7}.

\bigskip
Remind that Corollary \ref{CC1} immediately follows from Theorem~\ref{thm1},
 Corollary \ref{CC3} is a special  case of Theorem~\ref{thm3},
whereas Corollary \ref{C2+} simply follows from Theorem~\ref{thm2}.
All examples and lemmas are proved after their statements.
Thus, all results of the paper are proved.

\bigskip
 \noindent {\bf\large Acknowledgments} \quad
\noindent This work  was supported by the National Natural Science
Foundation of China (Grant No. 11931018), Tianjin Natural Science Foundation
and the program of fundamental scientific researches of the SB RAS № I.1.3., project № 0314-2019-0008.

\end{document}